\documentclass[11pt, reqno]{amsart}
 
 \usepackage{amsmath}
 \usepackage{amssymb}
 \usepackage{bm}

\DeclareMathAccent{\mathring}{\mathalpha}{operators}{"17}

 \usepackage{color}

\newcommand{\mysection}[1]{\section{#1}
      \setcounter{equation}{0}}

\newtheorem{theorem}{Theorem}[section]
\newtheorem{lemma}[theorem]{Lemma}

\theoremstyle{definition}
\newtheorem{assumption}{Assumption}[section]

\theoremstyle{remark}
\newtheorem{remark}{Remark}[section]

\newcommand\cbrk{\text{$]$\kern-.15em$]$}}
\newcommand\opar{\text{\raise.2ex\hbox{${\scriptstyle | }$}\kern-.34em$($} }
\newcommand{\tr}{\text{\rm tr}\,}

 \makeatletter
 \def\dashint{%
 \operatorname%
 {\,\,\text{\bf--}\kern-.98em\DOTSI\intop\ilimits@\!\!}}
\def\ninf{\qopname\relax\@empty{\!\!\phantom{p}\!\!\!inf}}
 \makeatother

\newcommand\bbeta{\text{\raise-.2ex\hbox{$\bm{\beta}$}}}

 \makeatletter
 \def\dashint{%
 \operatorname%
 {\,\,\text{\bf--}\kern-.98em\DOTSI\intop\ilimits@\!\!}}
\def\ninf{\qopname\relax\@empty{inf\phantom{p}\!\!\!}}
 \makeatother

\newcommand\bR{\mathbb{R}}

\newcommand\bS{\mathbb{S}}

\newcommand\cP{\mathcal{P}}

\newcommand\cH{\mathcal{H}}
\newcommand\cL{\mathcal{L}}

\newcommand\dist{{\rm dist}\,}

 \newcommand\esssup{\operatornamewithlimits{ess\,sup\,}}

\newcommand\supinf{\operatornamewithlimits{sup\,\,\,inf}}

\newcommand\supsup{\operatornamewithlimits{sup\,\,\,sup}}

\newcommand\osc{\operatornamewithlimits{osc\,}}

\begin{document}

\title[H\"older regularity of the first derivatives]
{On $C^{1+\alpha}$ regularity of solutions
of Isaacs parabolic equations
with VMO coefficients}

\author{N.V. Krylov}
\thanks{The  author was partially supported by
 NSF Grant DMS-1160569}
\email{krylov@math.umn.edu}
\address{127 Vincent Hall, University of Minnesota,
 Minneapolis, MN, 55455}



\begin{abstract}
We prove that  boundary value problems
for fully nonlinear second-order parabolic equations
admit  $L_{p}$-viscosity
solutions, which are in $C^{1+\alpha}$ for an $\alpha\in(0,1)$.
The equations have a special structure that the ``main'' part
containing only second-order derivatives is given by
a positive homogeneous function of second-order derivatives
and as a function of independent variables it is measurable
in the time variable and, so to speak, VMO in spatial
variables.
\end{abstract}

\maketitle

\mysection{Introduction}

In this article we take a function  $H(u,t,x)$,
$$
u=(u',u''),\quad
u'=(u'_{0},u'_{1},...,u'_{d}) \in\bR^{d+1},\quad u''\in\bS,\quad
(t,x)\in\bR^{d+1}, 
$$ 
where $\bS$ is the
set of symmetric $d\times d$ matrices, and
we are dealing with 
the parabolic  equation 
$$
\partial_t v(t,x)+H[v](t,x)
$$
\begin{equation}
                                                \label{7.29.1}
:=
 \partial_t v(t,x)+H(v(t, x),D v(t,x),D^{2}v(t, x),t, x)=f
\end{equation}
in subdomains of $(0,T)\times \bR^d $, where $T\in(0,\infty)$,
$$
\bR^{d}=\{x=(x_{1},...,x_{d}):x_{1},...,x_{d}\in \bR\},
$$
$$
\partial_t=\frac{\partial}{\partial t},\quad
 D^{2}u=(D_{ij}u),\quad Du=(D_{i}u),\quad
D_{i}=\frac{\partial}{\partial x_{i}},
\quad D_{ij}=D_{i}D_{j}.
$$

Our main goal is to establish the existence of $L_{p}$-viscosity
solutions of boundary value problems associated with \eqref{7.29.1},
solutions, which are in $C^{1+\alpha}$ for an $\alpha\in(0,1)$.
 
Let us briefly discuss what the author was able to find in the
literature concerning this kind of regularity. The articles
cited below contain a very large amount of information
concerning all kinds of issues in the theory
of fully nonlinear elliptic and parabolic equations,
but we will focus only on one of them.
Trudinger \cite{Tr88}, \cite{Tr89} and Caffarelli \cite{Ca89}
were the first authors who proved $C^{1+\alpha}$
regularity for fully nonlinear elliptic equations of type 
$$
F(u,Du,D^{2}u,x)
=f
$$ 
without convexity assumptions on $F$. The assumptions in these 
papers are different.
In \cite{Ca89} the function $F$ is independent of $u'$  
and, for each $u''$ uniformly sufficiently close to a function which is
 continuous
with respect to $x$. In \cite{Tr88} and \cite{Tr89} the function $F$
depends on all arguments but is H\"older continuous in $(u'_{0},x)$.
Next step in what concerns $C^{1+\alpha}$-estimates for the elliptic
case was done by    \'Swi{\c e}ch \cite{Sw97}, who considered
general $F$, imposed
the same condition as in \cite{Ca89} on the $x$-dependence, which is 
much weaker than in \cite{Tr88} and \cite{Tr89}, but also imposed the Lipschitz
condition on the dependence of $F$ on $u'_{1},...,u'_{d}$.
 In \cite{Tr88} and \cite{Tr89}
only continuity with respect to $u'_{1},...,u'_{d}$ is assumed.

In case of parabolic equations interior $C^{1+\alpha}$-regularity
was established by Wang \cite{Wa92} under the same kind of assumption
on the dependence of $H$ on $(t,x)$ as in \cite{Ca89} and assuming that
$H$ is almost independent of $u'_{1},...,u'_{d}$.
Then Crandall,  Kocan, and \'Swi{\c e}ch \cite{CKS00} generalized
the result of \cite{Wa92} to the case of full equation
again as in \cite{Sw97} assuming that 
$H$ is uniformly sufficiently close to a function which is continuous
with respect to $(t,x)$ and assuming the Lipschitz continuity
of $H$ with respect to $u'_{1},...,u'_{d}$
 and the continuity with respect to $u'_{0}$.

On the one hand, our class of equations is more narrow than the one
in \cite{CKS00} because we require the ``main'' part  of $H$, called $F$, be
positive homogeneous of degree one. On the other hand, we
do not require $H$ to be Lipschitz with respect to $u'_{1},...,u'_{d}$,
the continuity with respect to $u'$ suffices. Also
 we only need $F$ to be measurable
in $t$ and VMO in $x$, say, independent of $x$ and measurable
in $t$.

Our methods are absolutely different from the methods of above
cited articles. We do not use any 
ideas or facts from the theory of viscosity solutions.
Instead we rely on the methodology 
brought into the theory of fully nonlinear equations by Safonov
  \cite{Sa84},
 \cite{Sa88} and on an idea behind
 the proof of the main Lemma \ref{lemma 10.21.1} 
inspired by a probabilistic
interpretation of solutions of \eqref{7.29.1}.
We only focus on interior estimates of solutions
in smooth domains leaving to the interested reader
investigation of the same issues in nonsmooth domains
or  near the boundary of sufficiently regular ones. 

The article is organized as follows.
Section \ref{section 11.4.2} contains main results and 
some comments on them.
In Section \ref{section 11.4.3} we use a theorem from
 \cite{Kr13} to approximate the   equations with $H$
and with its main part $F$
by those for which the solvability is known.
We also leave to the interested reader carrying our
results over to elliptic equations.

In Section \ref{section 11.4.4} we show that
the approximate principal equation with $F$ admits solutions
locally well approximated in the sup norm by affine functions.
This is the most important part of the article.
Section \ref{section 11.4.5} contains estimates
of $C^{1+\alpha}$-norms of approximate equation with full $H$
and in Section \ref{section 11.4.6} we give the proof
of our main Theorem \ref{theorem 9.23.1}.
The last Section \ref{section 11.4.7} is actually an appendix,
which we need in order to be able to represent
positive homogeneous of order one functions depending
on parameters, such as $F$, as supinf's of affine functions whose
coefficients inherit the regularity properties
of the original function with respect to the parameters.

\mysection{Main results}
                                                       \label{section 11.4.2}

 To state our main results, we introduce
a few notation and assumptions.
Fix a constant
$\delta\in(0,1]$, and  set
$$
\bS_{\delta}=\{a\in\bS:\delta|\xi|^{2}\leq a_{ij}\xi_{i}\xi_{j}
\leq\delta^{-1}|\xi|^{2},\quad \forall\, \xi\in\bR^{d}\}, $$ where and
everywhere in the article the summation convention is enforced.

\begin{assumption}
                                    \label{assumption 9.23.1}
(i) The   function  $H(u,t,x)$
is measurable with respect to $(t,x)$ for any $u$ and
  Lipschitz continuous in $u''$ for every $u',(t,x)\in\bR^{d+1}$.

(ii) For any $(t,x)$, at all points of differentiability of 
$H(u,t,x)$  with
respect to $u''$, we have
 $
(H_{u''_{ij}}) \in \bS_{\delta}
 $.

(iii) There is a  function $\bar{H}(t,x)$ and a 
constant $K_{0}\geq0$ such that
$$
|H(u',0,t,x)|\leq K_{0}|u'|+\bar{H}(t,x).
$$

(iv) There is an increasing continuous function $\omega(r)$,
$r\geq0$, such that $\omega(0)=0$ and
$$
|H(u',u'',t,x)-H(v',u'',t,x)|\leq\omega(|u'-v'|)
$$
for all $u,v,t$, and $x$.

\end{assumption}
 
For $  R\in(0,\infty)$ and $(t,x)\in\bR^{d+1}$ introduce
$$
B_{R}=\{x\in\bR^{d}:|x|<R\},\quad B_{R}(x)=x+B_{R},
$$
$$
C_{ R}=(0,R^{2})\times B_{R},\quad
C_{R}(t,x)=(t,x)+C_{R}.
$$
For a Borel set $\Gamma$ in $\bR^{d+1}$  
by $|\Gamma|$ we  denote its Lebesgue measure.   Also
for a function $f$ on $\Gamma$ we set
$$
\dashint_{\Gamma}f(t,x)\,dxdt=\frac{1}{|\Gamma|} 
\int_{\Gamma}f(t,x)\,dxdt
$$
in case $\Gamma$ has a nonzero Lebesgue measure in $\bR^{d+1}$.
Similar notation is used in case of functions $f(x)$ on $\bR^{d}$.

We fix a constant $R_{0}\in(0,1]$ and 
for $\kappa\in(0,2]$ and measurable $f(t,x)$ introduce
$$
f_{\kappa}=\sup_{R\leq R_{0},t,x}
R^{2-\kappa}
\big(\dashint_{C_{R}(t,x)}|f(s,y)|^{d+1}\,dyds\big)^{1/(d+1)} .
$$
\begin{remark}
                                        \label{remark 10.23.1}
By H\"older's inequality for $p\ge d+1$
$$
\big(\dashint_{C_{R}(t,x)}|f(s,y)|^{d+1}\,dyds\big)^{1/(d+1)}
\leq NR^{-(d+2)/p}
\big(\int_{\bR^{d+1}}|f(s,y)|^{p}\,dyds\big)^{1/p},
$$
which shows that $f_{\kappa}<\infty$ if 
$f\in L_{p}(\bR^{d+1})$ and $\kappa\leq2-(d+2)/p$. It is
useful to observe that one can take $\kappa>1$ if
$f\in L_{p}(\bR^{d+1})$ for
  $p>d+2$.
\end{remark}
 
In the following assumption there are three
 objects $\kappa_{1}=\kappa(d,\delta)\in(1,2)$,
any $\kappa\in(1,\kappa_{1}]$, and 
 $\theta=\theta(\kappa,d,\delta)\in(0,1]$.
The values of $\kappa_{1}$ and $\theta$ are specified later
in the proof of Lemma \ref{lemma 10.23.3}.

 \begin{assumption}
                                  \label{assumption 10.23.1}
We have a representation 
$$
H(u,t,x)=F(u'',t,x)+G(u,t,x).
$$ 

(i) The functions $F$ and $G$ are measurable functions 
of their arguments.

(ii) For all values of the arguments
$$
|G(u,t,x)|\leq K_{0}|u'|+\bar{H}(t,x)
$$
\and there exists a $\kappa\in(1,\kappa_{1}]$ such that $\bar{H}_{\kappa}
<\infty$.

(iii) The function $F$ is positive homogeneous
of degree one with respect to $u''$,
is Lipschitz continuous with respect  to
$u''$, and at all points of differentiability of
$F$ with respect to $u''$ we have $F_{u''}\in\bS_{\delta}$.

(iv) For any $R\in(0,R_{0}]$, $(t,x)\in\bR^{d+1}$,
 and $u''\in\bS$
with $|u''|=1$ ($|u''|:=(\tr u''u'')^{1/2}$), we have
$$
\theta_{R,t,x}:=\dashint_{C_{R}(t,x)} 
 |F(u'',s,y)-\bar{F}_{R, x}(u'',s)|\,dsdy\leq\theta,
$$
where
$$
\bar{F}_{R, x}(u'',s)=\dashint_{B_{R}(x)}F(u'',s,y)\,dy.
$$

\end{assumption}

\begin{remark}
                                     \label{remark 10.25.1}

Assumption \ref{assumption 10.23.1} (ii) is stronger than
Assumption \ref{assumption 9.23.1} (iii) which is singled out
for methodological purposes.

Also observe that one can take $\theta=0$ 
in Assumption \ref{assumption 10.23.1} (iv)
if $F$ is 
independent of $x$.
\end{remark}

Fix a $T\in(0,\infty)$ and for domains $\Omega\in\bR^{d}$ define
$$
\Omega_{T}=(0,T)\times\Omega,\quad\partial'\Omega_{T}=\bar{\Omega}_{T}
\setminus(\{0\}\times\Omega).
$$
For $\kappa\in(0,1]$ and 
functions $\phi(t,x)$ on $\bar{\Omega}_{T}$ set
$$
[\phi]_{C^{ \kappa}(\bar{\Omega}_{T})}=
\sup_{(t,x),(s,y)\in
\bar{\Omega}_{T}}\frac{|\phi(t,x)-\phi(s,y)|}
{|t-s|^{\kappa/2}+|x-y|^{\kappa}},\quad
\|\phi\|_{C (\bar{\Omega}_{T})}=\sup_{\bar{\Omega}_{T}}|\phi|,
$$
$$
\|\phi\|_{C^{ \kappa}(\bar{\Omega}_{T})}=\|\phi\|_{C (\bar{\Omega}_{T})}
+[\phi]_{C^{ \kappa}(\bar{\Omega}_{T})}.
$$
For $\kappa\in(1,2]$ and sufficiently regular $\phi$ set
$$
[\phi]_{C^{ \kappa}( \bar{\Omega}_{T})}=
\sup_{t,s\in[0,t], x\in\bR^{d} }\frac{|\phi(t,x)-\phi(s,x)|}
{|t-s|^{\kappa/2} }
$$
$$
+\sup_{x,y\in\bar{\Omega} ,t\in[0,T]}\frac{|D\phi(t,x)-D\phi(t,y)|}
{ |x-y|^{\kappa-1}},\quad
\|\phi\|_{C^{ \kappa}( \bar{\Omega}_{T})}=
\|\phi\|_{C^{1}(\bar{\Omega}_{T})}+[\phi]_{C^{ \kappa}( \bar{\Omega}_{T})}.
$$
   The set of functions
with finite norm $\|\cdot\|_{C^{\kappa}(\bar{\Omega}_{T})}$ 
is denoted by $C^{\kappa}(\bar{\Omega}_{T})$.

\begin{remark}
                                   \label{remark 10.28.1}
According to the above notation $C^{2}(\bar{\Omega}_{T})$
is not what is usually meant. Therefore, we are going to use
the symbol
 $W^{1,2}_{\infty}(\Omega_{T})\cap C(\bar{\Omega}_{T})$ instead
for the space provided with norm
$\|\cdot\|_{C^{2}( \bar{\Omega}_{T})}$. One should keep
this in mind when we consider all $\kappa\in(0,2]$ at once.
\end{remark}

For   sufficiently regular functions $\phi(t,x)$ we set
\begin{equation}
                                              \label{10.31.1}
H[\phi](t,x)=H(\phi(t,x),D\phi(t,x),D^{2}\phi(t,x),t,x).
\end{equation}
 Similarly we introduce $F[\phi]$ and other operators
if we are given functions of $u,t,x$.

Everywhere below  $\Omega$ is a bounded $C^{2}$ domain in $\bR^{d}$
and $T\in(0,\infty)$.
The following is the main result of the paper. We refer the reader
to \cite{CKS00} for the definition of $L_{p}$-viscosity
solutions and their numerous properties.

\begin{theorem}
                                    \label{theorem 9.23.1}
Let  $g\in W^{1,2}_{\infty}(\Omega_{T})\cap C(\bar{\Omega}_{T})$.
Then  there is a function $v\in
 C^{\kappa}_{loc}(\Omega _T)\cap C(\bar{\Omega}_{T}) $
which, for any $p>d+2$, is an $L_{p}$-viscosity solution of 
the equation
\begin{equation}
                                               \label{9.23.2}
\partial_t v+ H[v] =0
\end{equation}
in $\Omega_T$ (a.e.) with boundary 
condition $v =g $ on $\partial'\Omega_{T}$.
 
Furthermore, for any $r,R\in(0,R_{0}]$ satisfying $r<R$ and $(t,x)
\in\Omega_{T}$ such that $C_{R}(t,x)\subset \Omega_{T}$ we have
\begin{equation}
                                                        \label{11.1.2}
[v]_{C^{\kappa}(C_{r}(t,x))}\leq N(R-r)^{-\kappa}\sup_{C_{R}(t,x)}|v|
+N \bar{H}_{\kappa} ,
\end{equation}
where $N$ depend  only on $d,\delta,K_{0}$, and $\kappa$
 (in particular, independent of $\omega$).

\end{theorem}

\begin{remark}
                                                      \label{remark 11.4.1}
A typical example of applications of Theorem \ref{theorem 9.23.1}
arises in connection with the theory of stochastic differential games 
where the so-called Isaacs equations play a major role.
To describe a particular case of these equations,
 assume that we are given countable sets
$A$ and $B$ and, for each $\alpha\in A$ and $\beta\in B$,
we have an
$\bS_{\delta}$-valued function $a^{\alpha\beta}(t,x)$
defined on $\bR^{d+1}$ and a
real-valued function  $G^{\alpha\beta}(u',t,x)$
defined for $u',(t,x)\in\bR^{d+1}$.
Suppose that these functions are measurable and
Assumption \ref{assumption 10.23.1} (ii)
is satisfied with $G^{\alpha\beta}$ in place of $G$
for any $\alpha\in A$ and $\beta\in B$ (and $\bar{H}$
independent of $\alpha\in A$ and $\beta\in B$).
Also suppose that Assumption \ref{assumption 9.23.1}
(iv) is satisfied with the same function $\omega$ and
with $G^{\alpha\beta}$ in place of $H$
for any $\alpha\in A$ and $\beta\in B$.
Finally, suppose that for any $R\in(0,R_{0}]$ and $(t,x)
\in\bR^{d+1}$
$$
\dashint_{C_{R}(t,x)} \supsup_{\alpha\in A\,\,\beta\in B}
|a^{\alpha\beta}(s,y)-\bar a^{\alpha\beta}(s)|\, dsdy\leq\theta,
$$
where
$$
\bar a^{\alpha\beta}(s)=\dashint_{B_{R} }a^{\alpha\beta}(s,y)\,dy.
$$

Upon introducing
$$
F(u'',t,x)=\supinf_{\alpha\in A\,\,\beta\in B}
 a^{\alpha\beta}_{ij}(t,x)u''_{ij},
$$
$$
G(u,t,x)=\supinf_{\alpha\in A\,\,\beta\in B}
\big[a^{\alpha\beta}_{ij}(t,x)u''_{ij}
+G^{\alpha\beta}(u',t,x)\big]-F(u'',t,x)
$$
one easily sees  that Theorem \ref{theorem 9.23.1} is applicable
to the equation
$$
\partial_{t}v+\supinf_{\alpha\in A\,\,\beta\in B}
\big[a^{\alpha\beta}_{ij}(t,x)D^{2}_{ij}v
+G^{\alpha\beta}(v,Dv, t,x)\big]=0.
$$

This example is close to the one from the introduction
in \cite{CKS00} and is more general, because $G^{\alpha\beta}$
are not assumed to be   linear in $u'$.
On the other hand, we suppose that 
Assumption \ref{assumption 10.23.1} (ii)
is satisfied with $G^{\alpha\beta}$ in place of $G$
uniformly in $\alpha,\beta$. In the situation of \cite{CKS00}
only
$$
\big(\supinf_{\alpha\in A\,\,\beta\in B}
G^{\alpha\beta}(0,\cdot,\cdot)\big)_{\kappa}<\infty
$$
is required.

\end{remark}

\begin{remark}
                                     \label{remark 10.25.2}

Assumption \ref{assumption 10.23.1} (iii), (iv)
can be replaced with the following which turns out
to be basically weaker (cf. \eqref{10.22.3}): There exist countable sets
$A$ and $B$ and functions $a^{\alpha\beta}(t,x)$
satisfying the conditions of Remark \ref{remark 11.4.1}
and there are numbers $f^{\alpha\beta}$ (independent of $(t,x)$) such that
$$
F(u'',t,x)=\supinf_{\alpha\in A\,\,\beta\in B}\big[
 a^{\alpha\beta}_{ij}(t,x)u''_{ij}+f^{\alpha\beta}\big]
\quad\text{and}\quad F(0,t,x)\equiv0.
$$
\end{remark}

\mysection{Auxiliary equations}
                                                       \label{section 11.4.3}

In the first result
of this section only Assumptions \ref{assumption 9.23.1}
is used.
By Theorem 2.1 of \cite{Kr13}
  there exists a convex positive homogeneous of degree
one function $P(u'')$ such that at all points of differentiability
of $P$ with respect to $u''$ we have $P_{u''}(u'')\in\bS_{\hat{\delta}}$,
where $\hat{\delta}=\hat{\delta}(d,\delta)\in(0,\delta/4)$
and such that the following fact holds in which 
 by $P[v]$ we mean a differential operator
constructed as in \eqref{10.31.1}.

\begin{theorem}
                                          \label{theorem 9.12.1}

Let  $K\geq0$ be a fixed constant,  
$g\in W^{1,2}_{\infty}(\Omega_{T})\cap C(\bar{\Omega}_{T})$.
 Assume that $\bar{H}$ is bounded. Then
the equation 
\begin{equation}
                                             \label{9.22.3}
\partial_{t}v+\max(H[v] ,P[v]-K)=0
\end{equation}
in $\Omega_{T}$  
with  boundary
condition $v=g$ on $\partial'\Omega_{ T}$ has a 
  solution $v\in
C(\bar{\Omega}_T)\cap W^{1,2}_{\infty,\text{loc}}(\Omega_{T})$.
In addition,  
$$
|v|, |D v|,\rho |D^{2} v |,|\partial_t v |\leq N(\sup_{\Omega_{T}}
\bar{H}+K
+\|g\|_{ C ^{1,2}(\Omega_T)})\quad\text{in} \quad \Omega_T \quad (a.e.),
$$
where
$$
\rho=\rho(x)=\dist(x,\bR^{d}\setminus\Omega) 
$$
 and  $N$ is a constant depending only on $\Omega$,
$T$, $K_{0}$, and $\delta$ (in particular, independent of $\omega$).

\end{theorem}

Theorem \ref{theorem 9.12.1} is applicable to
the equation
\begin{equation}
                                              \label{10.21.1}
\partial_{t}u+\max(F[u],P [u]-K)=0,
\end{equation}
which we want to rewrite in a different form.

First we
observe that if
in Section \ref{section 11.4.7} we take $B=\{0\}\times\bS_{\delta}$,
 take a strictly convex open set $B'_{0}$ in $\bS$ such that
$\bS_{\delta}\subset B'_{0}\subset \bS_{\delta/2}$, and set
$B_{0}=\{0\}\times \bar B'_{0}$, then by Theorem 
\ref{theorem 8.30.1} we have
$$
F(u'',t,x)=\supinf_{\alpha\in A_{1}\,\,\beta\in B}
a^{\alpha\beta}_{ij}(t,x)u''_{ij},\quad
\bar{F}(u'',t)=\supinf_{\alpha\in A_{1}\,\,\beta\in B}
\bar{a}^{\alpha\beta}_{ij}(t )u''_{ij}, 
$$
where $A_{1}=\bS$, for $\alpha\in A_{1}$ and 
$\beta=(0,\beta')\in B$,
$$
a^{\alpha\beta} (t,x)=\lambda^{\alpha\beta}(t,x)\beta'
+(1-\lambda^{\alpha\beta}(t,x))G_{u''}(\alpha)
$$
$$
\bar{a}^{\alpha\beta} (t)=\bar\lambda^{\alpha\beta}(t,x)\beta'
+(1-\bar\lambda^{\alpha\beta}(t,x))G_{u''}(\alpha),
$$
$$
G(u'')=\sup_{\beta'\in B'_{0}}\beta'_{ij}u''_{ij},
$$
$$
\lambda^{\alpha\beta}(t,x)=1\wedge
\frac{G(\alpha)-F(\alpha,t,x)}{G(\alpha)-
\beta'_{ij}\alpha_{ij}}\qquad\left(\frac{0}{0}=1\right),
$$
and $\bar\lambda^{\alpha\beta}(t)$ is defined similarly.
From Section \ref{section 11.4.7} we also know that,
for a constant $\mu>0$, we have $G(\alpha)-
\beta'_{ij}\alpha_{ij}\geq \mu|\alpha|$
if $\beta=(0,\beta')\in B$ and $\alpha\in A_{1}$.

Next, since $P(u'')$ is positive homogeneous, convex, and $P_{u''}\in\bS_{\hat
{\delta}}$, there exists a closed  set $A_{2}\subset \bS_{\hat
{\delta}}$ such that
$$
P(u'')=\sup_{\alpha\in A_{2}}\alpha_{ij}u''_{ij}.
$$
For uniformity of notation introduce $\hat{A}$ as a disjoint union
of $A_{1}$ and $A_{2}$ and for $\beta\in B$ and $\alpha\in A_{2}$
set
$$
a^{\alpha\beta}(t,x)=\bar a^{\alpha\beta}(t)=\alpha,\quad f^{\alpha\beta}=0.
$$

Also for $\alpha\in\hat{A}$ and $\beta\in B$ introduce
$\sigma^{\alpha\beta} (t,x)=[a^{\alpha\beta} (t,x)]^{1/2}$,
$\bar\sigma^{\alpha\beta} (t )=[\bar a^{\alpha\beta} (t )]^{1/2}$,
$$
L^{\alpha\beta}v(t,x)=a^{\alpha\beta}_{ij}(t,x)D_{ij}v(t,x),\quad
\bar L^{\alpha\beta}v(t,x)=\bar a^{\alpha\beta}_{ij}(t )D_{ij}v(t,x).
$$

Next we have the following which is essentially Remark 3.1
of \cite{DKL} with the proof based on the positive homogeneity and
 Lipschitz 
continuity of $F$ with respect to $u''$.

\begin{lemma}
                                         \label{lemma 10.6.4}
 There is a function $\theta=\theta(\mu)=\theta(\mu,d,\delta)>0$
defined for $\mu>0$ such that   Assumptions \ref{assumption 10.23.1}
(i), (iii), (iv) being satisfied with this $\theta(\mu)$ implies that
for any $R\in(0,R_{0}]$ and $(t,x)\in\bR^{d+1}$
$$
\dashint_{C_{R}(t,x)} \sup_{u''\ne 0}\frac{
 |F(u'',s,y)-\bar{F}_{R, x}(u'',s)|}{ |u''|}\,dsdy\leq\mu.
$$
\end{lemma}

Note that by Lemma \ref{lemma 10.6.4} and Theorem \ref{theorem 9.4.1}  
for any $R\in(0,R_{0}]$ and $(t,x)\in\bR^{d+1}$
\begin{equation}
                                                   \label{10.22.3}
\mu_{R,t,x}:=\dashint_{C_{R(t,x)}}\sup_{\alpha\in \hat A,\beta\in B}
  |a^{\alpha\beta}(s,y)-
\bar{a}^{\alpha\beta}(s) | \,dsdy\leq N\mu,
\end{equation}
where the constant $N$ depends only on $d$ and $\delta$.
On $\bS_{\hat\delta}$ the function $a^{1/2}$ is Lipschitz
continuous and therefore \eqref{10.22.3} also holds if we replace
$a$ with $\sigma$. 

Finally, observe that equation \eqref{10.21.1} is easily rewritten
as
\begin{equation}
                                               \label{10.28.2}
\partial_{t}u
+ \supinf_{\alpha\in\hat{A}\,\,\beta\in B}
\big[ L^{\alpha\beta}u( t,x)+f^{\alpha\beta}_{K} ]=0 ,
\end{equation}
where
$f^{\alpha\beta}_{K} = -KI_{\alpha\in A_{2}} $.

\mysection{Main estimate for solutions of \protect\eqref{10.21.1}}
                                                        \label{section 11.4.4}
Take $R,K\in(0,\infty)$ 
 and $g\in W^{1,2}_{\infty}(C_{R})\cap C(\bar C_{R})$. 
By Theorem \ref{theorem 9.12.1} there exists
  $u\in W^{1,2}_{\infty,loc}(C_{R})\cap C(\bar{C}_{R})$
such that $u=g$ on $\partial'C_{R}$ and equation \eqref{10.21.1}
holds (a.e.) in $C_{R}$.  By the maximum principle such $u$
is unique.

Here is the main result of this section.
\begin{theorem}
                                             \label{theorem 9.6.2}
There exist constants $\kappa_{0}
\in(1,2]$ and $N\in(0,\infty)$ depending only
on $d$ and $\delta$ such that for each $r\in(0,R]$
one can find an affine function $\hat{u}=\hat{u}(x)$
such that
$$
|u-\hat{u}|\leq  N  (\mu^{\kappa/(6d+6)} _{R}
\vee \mu^{1/(6d+6)} _{R})
[g]_{C^{ \kappa}(\bar C_{R})}R^{\kappa}+
Nr^{ \kappa_{0}}(R-r)^{-\kappa_{0}}\osc_{C_{R}}(g-\hat{g}) 
$$
  in $\bar C_{r}$ for any $\kappa\in(0,2]$,
where $\mu_{R}=\mu_{R,0,0}$ and
  $\hat{g}=\hat{g}(x)$ is any affine function of $x$.

\end{theorem}

By using parabolic dilations one easily sees that one may 
take  $R=1$.
In that case we first prove a few auxiliary results. 
Introduce $\bar{u}$ as a unique solution of
\eqref{10.21.1} (in $C_{1}$) with $\bar{F}$ in place of $F$ and the same
boundary condition on $\partial'C_{1}$.
Below by $N$ with occasional subscripts
 we denote various constants depending only
on $d$ and $\delta$.
\begin{lemma}
                                           \label{lemma 10.28.1}
Let $\kappa\in(0,2]$ and
\begin{equation}
                                               \label{10.28.5}
[g]_{C^{\kappa}(\bar C_{1})}=1.
\end{equation}
Then for any $\varepsilon>0$ there exists an
infinitely differentiable function 
$g^{\varepsilon}$ on $\bR^{d+1}$ such that in $\bar C_{1}$
\begin{equation}
                                               \label{10.28.3}
|g-g^{\varepsilon}|\leq N\varepsilon^{\kappa},
\quad
|\partial_{t}g^{\varepsilon}|+|D^{2}g^{\varepsilon}|
+\varepsilon|D^{3}g^{\varepsilon}|+
\varepsilon|D\partial_{t} g^{\varepsilon}|\leq
N\varepsilon^{\kappa-2},
\end{equation}
where $N$ depends only on $d$.
Furthermore, 
for $w=u,\bar u$   in $\bar C_{1}$ we have
\begin{equation}
                                               \label{10.28.4}
|w(t,x)-g^{\varepsilon}(t,x)|\leq N\varepsilon^{\kappa-2}
(1-|x|^{2})^{\kappa/2}+N\varepsilon^{\kappa} .
\end{equation}

\end{lemma}

Proof. The first assertion is well known and is obtained by first
continuing $g(t,x)$ as a function of $t$ to $\bR$ to become an even
2-periodic function, then continuing thus obtained function
across $|x|=1$ almost preserving \eqref{10.28.5} in the whole space
and then taking  convolutions
with $\delta$-like kernels.

Then, since $K\geq0$, for $w^{\varepsilon}
=u-g^{\varepsilon}$ we have
$$
\partial_{t}w^{\varepsilon}+
\partial_{t}g^{\varepsilon}+\max[F(D^{2}w^{\varepsilon}
+D^{2}g^{\varepsilon}),P(D^{2}w^{\varepsilon}
+D^{2}g^{\varepsilon})]\geq0,
$$
which in light of \eqref{10.28.3} implies that
$$
\partial_{t}w^{\varepsilon}+\max(F[w^{\varepsilon}],
  P[w^{\varepsilon}
 ])\geq-N_{1}\varepsilon^{\kappa-2}.
$$
Next, it is easily seen that there is a constant $N$
($=N(d,\delta 
 )$) such that for $\phi^{\varepsilon}(t,x)=
NN_{1}\varepsilon^{\kappa-2}(1-|x|^{2})$
we have
$$
\partial_{t}\phi^{\varepsilon}+\max(F[\phi^{\varepsilon}],
  P[\phi^{\varepsilon}
 ])\leq-N_{1}\varepsilon^{\kappa-2}
$$
in $C_{1}$. It follows
by the parabolic Alexandrov maximum principle that in $C_{1}$
\begin{equation}
                                              \label{10.28.6}
w^{\varepsilon}\leq \phi^{\varepsilon}
+\sup_{\partial'C_{1}}(w^{\varepsilon}- 
\phi^{\varepsilon}),\quad
u\leq g^{\varepsilon} +N \varepsilon^{\kappa-2}(1-|x|^{2})
+N\varepsilon^{\kappa},
\end{equation}
where $N$ depends only on $d$ and $\delta$. 

On the other hand, 
$$
\partial_{t}w^{\varepsilon}+
\partial_{t}g^{\varepsilon}+ F(D^{2}w^{\varepsilon}
+D^{2}g^{\varepsilon}) \leq0,
$$
$$
\partial_{t}w^{\varepsilon}+
  F[w^{\varepsilon}
 ]\leq N\varepsilon^{\kappa-2},
$$
and with perhaps different constant in the formula for
$\phi^{\varepsilon}$
$$
w^{\varepsilon}\geq-\phi^{\varepsilon} 
+\inf_{\partial'C_{1}}(\phi^{\varepsilon}+
 w^{\varepsilon}),\quad
u\geq g^{\varepsilon}   -N \varepsilon^{\kappa-2}(1-|x|^{2})
-N\varepsilon^{\kappa},
$$
which   along with
\eqref{10.28.6} yields \eqref{10.28.4} for $w=u$.
The proof of \eqref{10.28.4} for $w=\bar u$
is identical and the lemma is proved.

\begin{lemma}
                                           \label{lemma 10.21.1}
For
any $\kappa\in(0,2]$ in $\bar C_{1}$ we have
\begin{equation}
                                              \label{10.21.2}
|u-\bar{u}|\leq  N (\mu_{1}^{\kappa/(6d+6)}  \vee
\mu_{1}^{1/( d+1)})
[g]_{C^{ \kappa}(\bar C_{1})}.
\end{equation}
 
\end{lemma}

Proof. To simplify some formulas observe that if 
$[g]_{C^{\kappa}(\bar C_{1})}=0$, then $g$ is an affine function
of $x$ independent of $t$, so that $u=\bar u=g$ and 
we have nothing to prove. However, if $[g]_{C^{\kappa}(\bar C_{1})}>0$,
we can divide equation \eqref{10.21.1} by this quantity,
and, since our assertion means, in particular, that $N$
in \eqref{10.21.2} is independent of $K$, we can reduce the general 
situation to the one in which \eqref{10.28.5} holds.
 Therefore, below we assume \eqref{10.28.5}.

On sufficiently regular
 functions $u(t,x,\bar{x})$, $t\in\bR$, $ x,\bar x \in\bR^{d}$,
introduce
$$
\Phi[u](t,x,\bar{x})
=  \supinf_{\alpha\in\hat{A}\,\,\beta\in B}
\big[ \cL^{\alpha\beta}u(t,x,\bar{x})+f^{\alpha\beta}_{K} \big],
$$
where
$$
\cL^{\alpha\beta}u(t,x,\bar{x})=a^{\alpha\beta} _{ij}(t,x)
D^{x}_{ij}u(t,x,\bar{x}) 
+\hat{a}^{\alpha\beta} _{ij}(t,x,\bar x)D^{x\bar x}_{ij}u(t,x,\bar{x})
$$

$$
+\check{a}^{\alpha\beta} _{ij}(t,x,\bar x)D^{\bar x x}_{ij}u(t,x,\bar{x})
+\bar{a}^{\alpha\beta} _{ij}(t )D^{\bar x}_{ij}u(t,x,\bar{x}) ,
$$

$$
D^{x}_{ij}=\frac{\partial^{2}}{\partial x_{i}\partial x_{j}},\quad
D^{x\bar x}_{ij}=\frac{\partial^{2}}{\partial x_{i}\partial\bar x_{j}},\quad
D^{\bar x x}_{ij}=\frac{\partial^{2}}{\partial\bar x_{i}\partial x_{j}},\quad
D^{\bar x}_{ij}=
\frac{\partial^{2}}{\partial\bar x_{i}\partial\bar x_{j}},\quad
$$

$$
\hat{a}^{\alpha\beta} (t,x,\bar x)=
\sigma^{\alpha\beta}(t,x)\bar \sigma^{\alpha\beta}(t ),
\quad
\check{a}^{\alpha\beta} (t,x,\bar x)=
\bar \sigma^{\alpha\beta}(t )\sigma^{\alpha\beta}(t,x).
$$
 
Observe that for $\lambda,\bar\lambda\in\bR^{d}$ we have
$$
a^{\alpha\beta} _{ij}\lambda_{i}\lambda_{j}+
\hat{a}^{\alpha\beta} _{ij}\lambda_{i}\bar\lambda_{j}+
\check{a}^{\alpha\beta} _{ij}\bar\lambda_{i}\lambda_{j}
+\bar{a}^{\alpha\beta} _{ij}\bar\lambda_{i}\bar\lambda_{j}
=|\sigma^{\alpha\beta}\lambda+\bar\sigma^{\alpha\beta}\bar\lambda|^{2}\geq0,
$$
so that $\Phi$ is a (degenerate) elliptic operator.

Next let $w\in W^{1,2}_{d+1}(C_{1})\cap C(\bar C_{1})$ be a solution
of the equation
$$
\partial_{t}w+\supsup_{\alpha\in\hat A\,\,\beta\in B} L^{\alpha\beta}
w=-\supsup_{\alpha\in  A\,\,\beta\in B}|a^{\alpha\beta}-\bar{a}^{\alpha\beta}|
=:-h
$$
in $C_{1}$ with zero boundary condition on $\partial' C_{1}$.
Such a unique solution exists by Theorem 1.1 of \cite{DKL}  
and by the parabolic Alexandrov estimate and \eqref{10.22.3}
we have in $C_{1}$ that
\begin{equation}
                                                       \label{10.29.1}
0\leq w\leq N\mu_{1}^{1/(d+1)}.
\end{equation}
One of reasons we need the function $w$ is that, as is easy to see,
there is a $\lambda>0$ depending only on $d$ and $\delta$ such that
for all $\alpha,\beta$ on $C_{1}$ we have
$$
\partial_{t}(\lambda w(t,x)+|x-\bar x|^{2})+\cL^{\alpha\beta}(t,x,\bar x)
(\lambda w(t,x)+|x-\bar x|^{2})
$$

$$
=\lambda 
(\partial_{t}w+L^{\alpha\beta}w)(t,x)+2|\sigma^{\alpha\beta}(t,x)
-\sigma^{\alpha\beta}(t)|^{2}\leq 0,
$$
where the inequality follows from the fact that $a^{1/2}$ is
a Lipschitz continuous function on $\bS_{\hat\delta}$, so that
$|\sigma^{\alpha\beta}-\bar\sigma^{\alpha\beta}|^{2}\le N
|a^{\alpha\beta}-\bar a^{\alpha\beta}|^{2}\leq N
|a^{\alpha\beta}-\bar a^{\alpha\beta}|$.

After that we proceed in two steps.

{\em Step 1. Estimate of $u-\bar u$ from above}.
According to Lemma \ref{lemma 10.28.1}
  for $|\bar x|=1$ and $|x|\leq1$ we have
$$
u(t,x)\leq g^{\varepsilon}(t,x) +N \varepsilon^{\kappa-2}(1-|x|^{2})
+N\varepsilon^{\kappa} 
$$

$$
\leq g^{\varepsilon}(t,x) +N \varepsilon^{\kappa-2}
|x-\bar x|
+N\varepsilon^{\kappa},
$$
where
$$
\varepsilon^{\kappa-2}|x-\bar x|\leq
\varepsilon^{\kappa-4}|x-\bar x|^{2}+\varepsilon^{\kappa},
$$
so that
$$
u(t,x)\leq g^{\varepsilon}(t,x) +N \varepsilon^{\kappa-4}|x-\bar x|^{2}
+N\varepsilon^{\kappa}.
$$
This inequality  also obviously holds
if $|x|=1$, $|\bar x|\leq 1$ or if $t=1$ and $|x|,|\bar x|\leq 1$.
 This shows that for $\varepsilon\in(0,1)$
$$
u^{\varepsilon}(t,x,\bar x):=u(t,x)-[g^{\varepsilon}(t,x)-
g^{\varepsilon}(t,\bar x)
 +N 
\varepsilon^{\kappa-6}e^{1- t}|x-\bar x|^{2} 
+N\varepsilon^{\kappa}]\leq g^{\varepsilon}(t,\bar x)
$$
on $\partial'[(0,1)\times C_{1}^{2}]$. 
Actually, above we could have replaced $\varepsilon^{\kappa-6}$
with $\varepsilon^{\kappa-4}$ but later on we will need to deal with terms of order
$\varepsilon^{\kappa-6}|x-\bar x|^{2}$ anyway.
Also observe that for $\varepsilon\in(0,1)$,
$$
I^{\varepsilon}(t,x,\bar x):=\partial_{t}u^{\varepsilon}(t,x,\bar x)+\Phi[u^{\varepsilon}](t,x,\bar x)
=\partial_{t}u (t,x)+\partial_{t}g^{\varepsilon}(t,\bar x)-
\partial_{t}g^{\varepsilon}(t,x)
$$

$$
+N\varepsilon^{\kappa-6}e^{1- t}|x-\bar x|^{2}
+\supinf_{\alpha\in \hat A\,\,\beta\in B}\big[a^{\alpha\beta}_{ij}D_{ij}u(t,x)
+\bar{a}^{\alpha\beta}_{ij}D_{ij}g^{\varepsilon}(t,\bar x)
-a^{\alpha\beta}_{ij}D_{ij}g^{\varepsilon}(t,x)
$$
$$
-N\varepsilon^{\kappa-6}e^{1- t}|\sigma^{\alpha\beta}(t,x)-\bar{\sigma}^{\alpha\beta}
(t)|^{2}+f^{\alpha\beta}_{K}\big],
$$
where
$$
|\bar{a}^{\alpha\beta}_{ij}D_{ij}g^{\varepsilon}(t,\bar x)
-a^{\alpha\beta}_{ij}D_{ij}g^{\varepsilon}(t,x)|\leq 
|\bar{a}^{\alpha\beta}(t)-a^{\alpha\beta}(t,x)|\,|D^{2}g^{\varepsilon}(t,x)|
$$

$$
+N|D^{2}g^{\varepsilon}(t,x) -D^{2}g^{\varepsilon}(t,\bar x)|\leq N\varepsilon
^{\kappa-2}h(t,x)+N\varepsilon^{\kappa-3}|x-\bar x|,
$$

$$
|\partial_{t}g^{\varepsilon}(t,\bar x)-
\partial_{t}g^{\varepsilon}(t,x)|\leq N\varepsilon^{\kappa-3}|x-\bar x|,
$$

$$
- N\varepsilon^{\kappa-3}|x-\bar x|+N\varepsilon^{\kappa-6} |x-\bar x|^{2}
\geq -N\varepsilon^{\kappa}.
$$
It follows that
for $\varepsilon\in(0,1)$,
$$
I^{\varepsilon}(t,x,\bar x)\geq \partial_{t}u (t,x) 
+\supinf_{\alpha\in \hat A\,\,\beta\in B}
\big[L^{\alpha\beta} u(t,x)+f^{\alpha\beta}_{K}\big]
$$

$$
-N\varepsilon^{\kappa-6}h -N\varepsilon^{\kappa } =
-N_{1}\varepsilon^{\kappa-6}h -N_{1}\varepsilon^{\kappa }
$$
in $(0,1)\times C_{1}^{2}$.

On the other hand,
$$
\bar{u}(t,\bar x)\geq g^{\varepsilon}(t,\bar x)-N\varepsilon^{\kappa-2}(1-|\bar x|^{2})
-N\varepsilon^{\kappa},
$$
which implies that
$$
\bar u^{\varepsilon}(t,x,\bar x):=\bar{u}(t,\bar x)
+ N_{2}\varepsilon^{\kappa-6}((1+\lambda  ) w(t,x)+ |x- \bar x|^{2})
$$

$$
+N_{2}(2-t)\varepsilon^{\kappa}\geq g^{\varepsilon}(t,\bar x)
$$
on $\partial'[(0,1)\times C_{1}^{2}]$. It is also easily seen that
by increasing $N_{2}$ if needed (which does not violate the above inequality)
we may assume that
in $(0,1)\times C_{1}^{2}$
$$
\partial_{t}\bar u^{\varepsilon}(t,x,\bar x)
+\Phi[\bar u^{\varepsilon}](t,x,\bar x)\leq  
-N_{1}\varepsilon^{\kappa-6}h -N_{1}\varepsilon^{\kappa }.
$$
Hence, by the maximum principle (see, for instance, Theorem 2.1 of
\cite{Kr78} or Theorem 3.4.2 of \cite{Kr87}) in 
$[0,1]\times\bar C_{1}^{2}$ we have
$$
\bar{u}(t,\bar x)
+ N\varepsilon^{\kappa-6} ( w(t,x)+|x-\bar x|^{2})
+N\varepsilon^{\kappa} \geq u(t,x)
$$

$$
-[g^{\varepsilon}(t,x)-
g^{\varepsilon}(t,\bar x) 
 +N 
\varepsilon^{\kappa-6}|x-\bar x|^{2} 
+N\varepsilon^{\kappa}]  ,
$$
which for $x=\bar x$ in light of \eqref{10.29.1} yields
$$
u(t,x)-\bar{u}(t,\bar x)\leq N(\varepsilon^{\kappa}+ 
\varepsilon^{\kappa-6}\mu_{1}^{1/(d+1)}).
$$
If $\mu_{1}\leq 1$, then for
 $\varepsilon=\mu_{1}^{1/(6d+6)}$ ($\leq1$) we get
 $u-\bar u\leq N\mu_{1}^{\kappa/(6d+6)}  $ and if 
$\mu_{1}\geq 1$, then for $\varepsilon=1$ we obtain
$u-\bar u\leq N\mu_{1}^{1/( d+1)}$, so that generally
$$
u-\bar u\leq N(\mu_{1}^{\kappa/(6d+6)}\vee\mu_{1}^{1/( d+1)}).
$$

{\em Step 2. Estimate of $u-\bar u$ from below}. Notice that
$$
\bar v^{\varepsilon}(t,x,\bar x):=
\bar{u}(t,\bar x)
- N\varepsilon^{\kappa-4}(\lambda w(t,x)+|x- \bar x|^{2})
-N\varepsilon^{\kappa}\leq g^{\varepsilon}(t,\bar x)
$$
on $\partial'[(0,1)\times C_{1}^{2}]$. It is also easily seen that
in $(0,1)\times C_{1}^{2}$
$$
\partial_{t}\bar v^{\varepsilon}(t,x,\bar x)
+\Phi[\bar v^{\varepsilon}](t,x,\bar x)\geq0.
$$
 
On the other hand, 
$$
v^{\varepsilon}(t,x,\bar x):=u(t,x)-[g^{\varepsilon}(t,x)-
g^{\varepsilon}(t,\bar x)]
$$
$$
 +N 
\varepsilon^{\kappa-6}e^{1-t}((1+\lambda) w(t,x)+|x-\bar x|^{2} )
+N(2-t)\varepsilon^{\kappa} \geq g^{\varepsilon}(t,\bar x)
$$
on $\partial'[(0,1)\times C_{1}^{2}]$ and the above computations
show that (for sufficiently large $N$)
$$
\partial_{t}  v^{\varepsilon}(t,x,\bar x)
+\Phi[  v^{\varepsilon}](t,x,\bar x)\leq0
$$
in $(0,1)\times C_{1}^{2}$. By the maximum principle $\bar v^{\varepsilon}
\leq v^{\varepsilon} $, which leads to the desired estimate of 
$u-\bar u$ from below and the lemma is proved.

\begin{lemma}
                                            \label{lemma 8.29.1} 
There exist constants $\kappa_{0}\in(1,2]$ and $N\in(0,\infty)$
depending only on $d$ and $\delta$
such that for any $r\in(0,1)$  
\begin{equation}
                                                   \label{8.29.1}
[\bar{u}]_{C^{ \kappa_{0} }(\bar C_{r})}
\leq N(1-r)^{-\kappa_{0}}\osc_{\bar C_{1 }}(g-\hat{g}),
\end{equation}
where $\hat{g}=\hat{g}(x)$ is any affine function of $x$.

\end{lemma}

Proof. First observe that $\bar{u}-\hat{g}$ satisfies the same
equation as $\bar{u}$ and the $C^{ \kappa }(\bar C_{r})$-seminorms
of these functions coincide if $\kappa\in(1,2]$.
 It follows that we may concentrate
on $\hat{g}\equiv0$.

For any $ \rho\in(0,1) $ the function
$\delta_{h}\bar{u}$ satisfies an equation of type
$$
\partial_{t}\delta_{l,h}\bar{u}+a_{ij}D_{ij}\delta_{l,h}\bar{u}=0
$$
in $C_{\rho}$
with some $(a_{ij})$ taking values in $\bS_{\hat{\delta}}$
 if $h$ is sufficiently small. By Corollary
4.3.6 of \cite{Kr87} for such $h$ and $r\in(0,\rho)$ we have
$$
[\delta_{l,h}\bar{u}]_{C^{ \gamma}(\bar C_{ r})}
\leq N(\rho-r)^{-\gamma}\sup_{\bar C_{\rho}}|\delta_{l,h}\bar{u}|,
$$
where $N$ and $\gamma\in(0,1)$ depend only on $\delta$
and $d$. By letting $h\to0$ we conclude
\begin{equation}
                                               \label{9.4.1}
[D\bar{u}]_{C^{ \gamma}(\bar C_{ r})}
\leq N(\rho-r)^{-\gamma}\sup_{\bar C_{\rho}}|D\bar{u}|.
\end{equation}

Next observe that for any function $f(x)$ of one variable
$x\in[0,\varepsilon]$, $\varepsilon>0$,  we have
$$
|f'(0)|\leq |f'(0)-(f(\varepsilon)-f(0)/\varepsilon|+
\varepsilon^{-1}\osc_{[0,\varepsilon)}f
\leq\varepsilon^{\gamma}[f']_{C^{\gamma}[0,\varepsilon]}+
\varepsilon^{-1}\osc_{[0,\varepsilon]}f.
$$
By applying this fact to functions $v(x)$ given in $B_{1}$ 
we obtain that for any $r_{n+1}<r_{n+2} \leq1$
and  any $\varepsilon\in(0,1)$
\begin{equation}
                                              \label{11.1.1}
|Dv|\leq\varepsilon^{\gamma}(r_{n+2}-r_{n+1 })^{\gamma} 
[Dv]_{C^{\gamma}(\bar B_{r_{n+2}})}
+\varepsilon^{-1}(r_{n+2}-r_{n+1})^{-1}
 \osc_{\bar B_{1}}v
\end{equation}
in $\bar B_{r_{n+1}}$.
 
Coming back to \eqref{9.4.1} and setting
$$
r_{0}=r,\quad
r_{n}=r+(1-r)\sum_{k=1}^{n}2^{-k},\quad n\geq1,
$$
we conclude
$$
A_{n}:=\sup_{[0,r_{n}^{2}]}[D\bar{u}(t,\cdot)
]_{C^{ \gamma}(\bar B_{ r_{n}})}
\leq N(r_{n+1}-r_{n})^{-\gamma}\sup_{\bar C_{ r_{n+1}}}|D\bar{u}|
$$
\begin{equation}
                                               \label{9.4.2}
\leq N_{1}\varepsilon^{\gamma} 
A_{n+2}+N_{2}(1-r)^{-(1+\gamma)}\varepsilon^{-1}2^{(1+\gamma)n}
\osc_{\bar C_{ 1}}\bar{u},
\end{equation}
where the constants $N_{i}$ are different from the one in \eqref{9.4.1}
but still depend  only on $\delta$ and $d$.
We first take $\varepsilon$ so that 
$$
N_{1}\varepsilon^{\gamma} 
=2^{-5},
$$
then take $n=2k$, $k=0,1,...$,  multiply both parts
of \eqref{9.4.2} by $2^{-5k}$ and sum up with respect to $k$.
Then upon observing that $(1+\gamma)2k\leq 4k$ we get
$$
\sum_{k=0}^{\infty}A_{2k}2^{-5 
k}\leq
\sum_{k=1}^{\infty}A_{2k}2^{-5k}+N(1-r)^{-(1+\gamma)}
\sum_{k=0}^{\infty}
2^{-k}\osc_{\bar C_{ 1}}\bar{u}.
$$
By canceling (finite) like terms we find
\begin{equation}
                                               \label{9.4.3}
\sup_{[0,r^{2}]}[D\bar{u}(t,\cdot)]_{C^{ \gamma}(\bar B_{r})}
\leq N(1-r)^{-(1+\gamma)}\osc_{\bar C_{ 1}}\bar{u}.
\end{equation}

Next, we use the fact that $\bar{u}$ itself satisfies the equation
$$
0=\partial_{t}\bar{u}+\max(\bar{F}[\bar{u}],\bar{P}[\bar{u}]-K)
-\max(0,-K)=\partial_{t}\bar{u}+a_{ij}D_{ij}\bar{u} 
$$
with some $(a_{ij})$ taking values in $\bS_{\hat{\delta}}$.
 Furthermore, for any $T\in(0,r^{2}]$
and $|x_{0}|\leq r$
the function
$v(t,x):=\bar{u}(t,x)-(x_{i}-x_{0i})D_{i}\bar{u}(T,x_{0})$
satisfies the same equation and
$$
|v(T,x)-v(T,x_{0})|\leq [D\bar{u}(T,\cdot)]_{C^{\gamma}
(\bar B_{\rho})}
|x-x_{0}|^{1+\gamma}
$$
for $|x-x_{0}|\leq \rho-r$, where $\rho=(1+r)/2$. 
Therefore, by Lemma 4.4.2 of \cite{Kr87},
applied with $R=\rho-r=(1-r)/2$ there,
for $t\in[0,T]$ we have
$$
|\bar{u}(t,x_{0})-\bar{u}(T,x_{0})|\leq
N[D\bar{u}(T,\cdot)]_{C^{\gamma}
(\bar B_{\rho})}(T-t)^{(1+\gamma)/2}
$$
$$
\leq N
 (1-r)^{-(1+\gamma)}
(T-t)^{(1+\gamma)/2}\osc_{\bar C_{ 1}}\bar{u}.
$$
This provides the necessary estimate of the oscillation of $\bar{u}$
in the time variable and along with \eqref{9.4.3} shows that
$$
[\bar{u}]_{C^{ 1+\gamma}(\bar C_{r})}
\leq N(1-r)^{-(1+\gamma)}\osc_{\bar C_{  1}}\bar{u}.
$$

Now the assertion of the lemma follows
from the fact that 
$$
\osc_{\bar C_{ 1}}\bar{u}=\osc_{\bar C_{ 1}}g.
$$
The lemma is proved.

{\bf Proof of Theorem \ref{theorem 9.6.2}}. Take $\hat{u}(t,x)
=\bar{u}(0,0)+x_{i}D_{i}\bar{u}(0,0)$ and
observe that in $C_{r}$
$$
|u-\hat{u}|\leq |u-\bar{u}|+|\bar{u}-\hat{u}|
\leq N (\mu_{1} ^{\kappa/(6d+6)}\vee \mu_{1} ^{1/(6d+6)})  
[g]_{C^{ \kappa}(\bar C_{1})}+I,
$$
where 
$$
I=|\bar{u}-\hat{u}|\leq 2r^{\kappa_{0}}[\bar{u}]_{C^{\kappa_{0}}(\bar C_{r})}
\leq Nr^{\kappa_{0}}(1-r)^{-\kappa_{0}}\osc_{\bar C_{1}}(g-\hat{g})
$$
so that the theorem is proved.

\mysection{Estimating $C^{\kappa}$-norm of solutions
of \protect\eqref{9.22.3}}
                                                       \label{section 11.4.5}

In this section we assume that $\bar{H}$ is bounded and investigate
solutions of \eqref{9.22.3} which exist by Theorem \ref{theorem 9.12.1}.
We take $\kappa_{0}\in(1,2]$ from Theorem \ref{theorem 9.6.2},
take a $\mu\in(0,1]$, and suppose that Assumption \ref{assumption 10.23.1}
(iv) is satisfied with $\theta=\theta(\mu)$ so that
\eqref{10.22.3} holds 
for any $R\in(0,R_{0}]$ and $(t,x)\in\bR^{d+1}$.

\begin{lemma}
                                                \label{lemma 10.22.3}
Let $R\in(0,R_{0}]$  and let
$ v \in W^{1,2}_{\infty}(\bar C_{R})\cap C(\bar{C}_{R})$ be a   solution
of \eqref{9.22.3} in $\bar C_{R}$.
Then for each $r\in(0,R)$
  one can find
 an affine function $\hat{v}(x)$
such that in $C_{r}$ for any $\kappa\in[1,2]$
$$
|v-\hat{v}|\leq  
N \mu^{\kappa/(6d+6)} 
[v]_{C^{ \kappa}(\bar C_{R})}R^{\kappa}+
Nr^{ \kappa_{0}}(R-r)^{-\kappa_{0}}R^{\kappa}[v]_{C^{ \kappa}(\bar C_{R})}
$$
$$
 +
NK_{0}R^{2}\sup_{\bar C_{R}}(|v|+|Dv|)+NR^{\kappa}  
 \bar{H}_{\kappa} ,  
$$
where   the
constants $N$ depend  only on $d$ and $\delta$.

\end{lemma}

Proof. 
Observe that
$$
\max(H[v] ,P[v]-K)=\max(F[v] ,P[v]-K)+h
$$
where $h$ defined by the above equality satisfies
$$
|h|\leq|H[v]-F[v]| \leq K_{0}(|v|+|Dv|)+\bar{H} .
$$
Next define $ u \in W^{1,2}_{d+1}(C_{R})
\cap C(\bar{C}_{R})$ as a unique solution
$$
\partial_{t}u+\max(F[u],P[u]-K)=0
$$
with  boundary data $u=v$ on $\partial' C_{R}$.
Then  there exists an 
$\bS_{\hat{\delta}}$-valued function $a$ such that
in $C_{R}$ we have
$$
\partial_{t}(v-u)+a_{ij}D_{ij}(v-u)+h=0.
$$
By the parabolic Alexandrov estimate
(cf. our comment concerning this estimate in
a more general situation in the proof of Lemma \ref{lemma 9.20.1})
$$
|v-u|\leq NR^{d/{d+1}}\|h\|_{L_{d+1}(C_{R})}=NR^{2}
\big(\dashint_{C_{R}}|h|^{d+1}\,dxdt\big)^{1/(d+1)} 
$$
$$
\leq NK_{0}R^{2}\sup_{\bar C_{R}}(|v|+|Dv|)+NR^{\kappa} \bar{H}_{\kappa}.
$$

After that our assertion follows from 
  Theorem \ref{theorem 9.6.2} and the lemma is proved.  

 Here is a result, which can be easily extracted
from the proof  of Theorem 2.1 of \cite{Sa88}.

\begin{lemma}
                                               \label{lemma 10.23.1}
Let $r_{0}\in(0,\infty)$, $\kappa\in(1,2)$, $\phi\in 
C^{\kappa}(\bar C_{r_{0}})$ and assume that there is a constant
$N_{0}$ such that for any $(t,x)\in C_{r_{0}}$ 
and $r\in(0,2r_{0}]$
there exists an affine function $\hat{\phi}=\hat{\phi}(x)$ such that
$$
\sup_{\bar C_{r}(t,x)\cap
\bar C_{r_{0}}} |\phi-\hat{\phi}|\leq N_{0}r^{\kappa}.
$$
Then
$$
[\phi]_{C^{\kappa}(\bar C_{r_{0}})}\leq NN_{0},
$$
where $N  $ depends only on $d$ and $\kappa$.

\end{lemma}

\begin{lemma} 
                                          \label{lemma 10.23.3}
Take $r_{1}\in(0,R_{0}]$, $r_{0}\in(0,r_{1})$,
and define
$$
\kappa_{1} =
\frac{1+\kappa_{0}}{2}.
$$
Let 
$ v \in W^{1,2}_{\infty}(C_{r_{1}})\cap C(\bar{C}_{r_{1}})$ be a 
  solution
of \eqref{9.22.3} in $C_{r_{1}}$
and let $\kappa\in(1,\kappa_{1}]$.
 Then there exists $\theta=\theta(\kappa,d,\delta)\in(0,1]$ such that,
if Assumption \ref{assumption 10.23.1} (iv) is satisfied with this
$\theta$, then  
$$
[v]_{C^{\kappa}(\bar C_{r_{0}})}
\leq (1/2)[v]_{C^{  \kappa}(\bar C_{r_{1}})}
+N(K_{0}+1) (r_{1}-r_{0})^{-\kappa }\sup_{\bar C_{r_{1}}} |v|
$$
\begin{equation}
                                                   \label{10.23.4}+
N (K_{0}+1)(r_{1}-r_{0})^{-(\kappa-1) }\sup_{\bar C_{r_{1}}} |Dv| 
+  N \bar{H}_{\kappa }  ,
\end{equation}
where $N =N (d,\delta,\kappa)$. 
\end{lemma}

Proof. To specify $\theta$ we first  
take a $\mu\in(0,1]$ and suppose that Assumption \ref{assumption 10.23.1}
(iv) is satisfied with $\theta=\theta(\mu)$ so that
\eqref{10.22.3} holds 
for any $R\in(0,R_{0}]$ and $(t,x)\in\bR^{d+1}$.

Then take $(t_{0},x_{0})\in C_{r_{0}}$, $\varepsilon\in(0,1)$, define 
$$
r'_{0}=\frac{\varepsilon}{3}(r_{1}-r_{0}), 
$$
and notice that for any $(t,x)\in C_{r_{0}'}(t_{0},x_{0})$,
  $r\in (0,2r_{0}']$, and $R=\varepsilon^{-1}r$,   we have
$$
C_{R}(t,x)\subset C_{r _{1}}.
$$
  Therefore,
by Lemma \ref{lemma 10.22.3} we can find an affine function
$\hat{v}(x)$ such that
$$
\sup_{C_{r}(t,x)\cap C_{r'_{0}}(t_{0},x_{0})}
|v-\hat{v}|\leq
\sup_{\bar C_{r}(t,x) }
|v-\hat{v}|
$$
$$
\leq N \mu^{\kappa /(6d+6)} 
[v]_{C^{  \kappa}(\bar C_{R}(t,x))} \varepsilon ^{-\kappa }
r^{\kappa }+
N\varepsilon^{ \kappa_{0}-\kappa }(1-\varepsilon)^{-\kappa_{0}}
r^{\kappa }[v]_{C^{  \kappa}(\bar C_{R}(t,x))}
$$

$$
 +
NK_{0}\varepsilon^{-2}
r^{2}\sup_{\bar C_{R}(t,x)}(|v|+|Dv|)+N
\varepsilon^{-\kappa }
r^{\kappa }\bar{H}_{\kappa }\leq N r^{\kappa }I(\varepsilon,r_{1}),  
$$
where the constants $N$ depend only on $d$ and $\delta$ and
$$
I(\varepsilon,r_{1}):=\big(\mu^{\kappa /(6d+6)}
\varepsilon ^{-\kappa }+\varepsilon^{ \kappa_{0}-\kappa }
(1-\varepsilon)^{-\kappa_{0}}
\big)[v]_{C^{  \kappa}(\bar C_{r_{1}})}  
$$

$$
+\varepsilon^{-2}K_{0}\sup_{\bar C_{r_{1}}}(|v|+|Dv|)
+\varepsilon^{-\kappa }
 \bar{H}_{\kappa }.
$$
It follows by Lemma \ref{lemma 10.23.1} that
$$
[v]_{C^{\kappa }(\bar C_{r'_{0}}(t_{0},x_{0}))}\leq N_{1}I(\varepsilon,r_{1}),
$$
where   $N_{1}$ depends only on $d$, $\kappa$, and $\delta$. We can now specify 
$\theta$ and $\varepsilon$. First we chose
$\varepsilon\in(0,1)$ so that
$$
N_{1}\varepsilon^{ \kappa_{0}-\kappa }
(1-\varepsilon)^{-\kappa_{0}}= 1/4.
$$
Since $\kappa_{0}-\kappa\geq (\kappa_{0}-1)/2>0$ and $\kappa_{0}$
depends only on $d$ and $\delta$ and $N_{1}$
depends only on $d$, $\kappa$, and $\delta$,
 $\varepsilon$ also 
depends only on $d$, $\kappa$, and $\delta$.
After that we take   $\mu=\mu(d,\kappa,\delta)\in(0,1]$ so that
$$
N_{1}\mu^{1/(6d+6)}
\varepsilon ^{-2}\leq 1/4
$$
and set $\theta=\theta(\mu(d,\kappa,\delta))$.
Then
\begin{equation}
                                              \label{10.25.1}
[v]_{C^{\kappa }(\bar C_{r'_{0}}(t_{0},x_{0}))}\leq 
(1/2)[v]_{C^{  \kappa}(\bar C_{r_{1}})}
+NJ,
\end{equation}
where $N=N(d,\delta,\kappa)$ and
$$
J=K_{0}\sup_{\bar C_{r_{1}}}(|v|+|Dv|)
+  \bar{H}_{\kappa }.
$$

Now observe that if $(t,x),(s,x)\in C_{r_{0}}$ and $t> s$, then
 either $|t-s|\leq (r_{0}')^{2}$, in which case 
$(t,x)\in C_{r'_{0}}(s,x)$ and
$$
(t-s)^{-\kappa /2}|v(t,x)-v(s,x)|
\leq (1/2)[v]_{C^{  \kappa}(\bar C_{r_{1}})}
+NJ
$$
owing to \eqref{10.25.1},
 or
  $|t-s|\geq (r_{0}')^{2}$ when 
$$
|v(t,x)-v(s,x)|\leq 2(t-s)^{\kappa /2}(r'_{0})^{- \kappa }
\sup_{\bar C_{r_{1}}}|v|\leq N(t-s)^{\kappa /2}(r_{1}-r_{0})
^{- \kappa }\sup_{\bar C_{r_{1}}}|v|.
$$

Next if $(t,x),(t,y)\in C_{r_{0}}$ and $x\ne y$, then either
$|x-y|<r_{0}'$, in which case $(t,y)\in  C_{r'_{0}}(t,x)$ and
$$
|x-y|^{-(\kappa-1) }|Dv(t,x)-Dv(t,y)|\leq
(1/2)[v]_{C^{  \kappa}(\bar C_{r_{1}})}+NJ,
$$
or else $|x-y|\geq r_{0}'$ and
$$
|Dv(t,x)-Dv(t,y)|\leq 2|x-y|^{\kappa-1 }(r_{0}')^{-(\kappa-1) }
\sup_{\bar C_{r_{1}}}|Dv|
$$
$$
\leq N|x-y|^{\kappa-1 }
(r_{1}-r_{0})
^{-(\kappa-1) }\sup_{\bar C_{r_{1}}}|Dv|.
$$
This proves \eqref{10.23.4} and the lemma.

\begin{theorem}
                                      \label{theorem 11.1.1}
Take $0<r<R\leq R_{0}$ and take $\kappa_{1}$,
$\kappa\in(1,\kappa_{1}]$, and $\theta$
from Lemma \ref{lemma 10.23.3}.
Let 
$ v \in W^{1,2}_{\infty}(C_{R})\cap C(\bar{C}_{R})$ be a 
  solution
of \eqref{9.22.3} in $C_{R}$. Then 
\begin{equation}
                                                        \label{11.1.02}
[v]_{C^{\kappa}(\bar C_{r})}\leq N(R-r)^{-\kappa}\sup_{\bar C_{R}}|v|
+N \bar{H}_{\kappa} ,
\end{equation}
where $N$ depends only on $d,\delta,K_{0}$, and $\kappa$.

\end{theorem}

Proof.
We proceed as in the proof of Lemma \ref{lemma 8.29.1}.
Fix a number $c\in(0,1)$ such that $c^{4}>3/4$ and   introduce
$$
r_{0}=r,\quad
r_{n}=r+c_{0}(R-r)\sum_{k=1}^{n}c^{k},\quad n\geq1,
$$
where $c_{0}$ is chosen in such a way that $r_{n}\to R$ as $n\to\infty$.
Then 
Lemma \ref{lemma 10.23.3} and \eqref{11.1.1} allow us to find constants
$N_{1}$ and $N$ depending only on $d,\delta,K_{0}$,
 and $\kappa$ such that for all $n$
and $\varepsilon\in(0,1)$
$$
A_{n}:=[v]_{C^{\kappa}(\bar C_{r_{n}})}\leq
(2^{-1}+N_{1}\varepsilon^{\kappa-1})A_{n+2}
$$
$$
+
N(R-r)^{-\kappa}
c^{-n\kappa}(1+\varepsilon^{-1})\sup_{\bar C_{R}}|v|+N\bar{H}_{\kappa}.
$$
we choose $\varepsilon< 1$ so  
that $2^{-1}+N_{1}\varepsilon^{\kappa-1}
\leq3/4$ and then recalling that $\kappa\leq 2$ conclude that
$$
\sum_{k=0}^{\infty}(3/4)^{k}A_{2k}\leq
\sum_{k=1}^{\infty}(3/4)^{k}A_{2k} +N\bar{H}_{\kappa}
$$
$$
+N(R-r)^{-\kappa}\sup_{\bar C_{R}}|v|
\sum_{k=0}^{\infty}(3/4)^{k} c^{-4k},
$$
where the last series converges since $3c^{-4}/4<1$.
By canceling like terms we come to \eqref{11.1.02}
and the theorem is proved.

\mysection{Proof of Theorem \protect\ref{theorem 9.23.1}}
                                                       \label{section 11.4.6}
First assume that $\bar{H}$ is bounded.
For $K>0$ denote by $v_{K}$ the solution of
 \eqref{9.22.3} with boundary condition $v=g$
on $\partial'\Omega_{T}$. By Theorem \ref{theorem 9.12.1}
such a solution exists is continuous in $\bar\Omega_{T}$
and has locally bounded derivatives.

Then the beginning of the proof of Lemma \ref{lemma 10.22.3}
shows that for an $\bS_{\hat{\delta}}$-valued function $(a_{ij})$
we have
$$
|\partial_{t}v_{K}+a_{ij}D_{ij}v_{K}|\leq K_{0}
(|v_{K}|+|Dv_{K}|+\bar{H} ),
$$
and the parabolic Alexandrov estimate shows that
\begin{equation}
                                                         \label{11.1.4}
|v_{K}|\leq N(\|g\|_{C (\Omega_{T})}+\|\bar{H}\|_{L_{d+1}(\Omega_{T})}
 ),
\end{equation}
where $N$ depends only on  $d$, $\delta$, $K_{0}$, and the diameter
of $\Omega$. 

Also
$$
|\partial_{t}(v_{K}-g)+a_{ij}D_{ij}(v_{K}-g)|
$$
\begin{equation}
                                                         \label{11.1.5}
\leq K_{0}
(|v_{K}|+|D(v_{K}-g)|)+\bar{H} +N (|\partial_{t}g|+|D^{2}g|+|Dg|),
\end{equation}
which, after we continue $(v-g)(t,x)$ for $t\geq T$ as zero,
by Theorem 4.2.6 of \cite{Kr87} yields that there exists an $\alpha=\alpha
(d,\delta)\in(0,1)$ such that for
any  domain $\Omega'\subset\bar\Omega'\subset\Omega$
\begin{equation}
                                                         \label{11.1.6}
|v_{K}|_{C^{\alpha}(\Omega'_{T})}\leq N,
\end{equation}
where $N$ depends only on the distance between the boundaries of
$\Omega'$ and $\Omega$ and on $T$, $d$, $\delta$, $K_{0}$,  the diameter
of $\Omega$, and the $L_{d+1}(\Omega_{T})$-norms of $\bar{H}$ and 
$|\partial_{t}g|+|D^{2}g|+|Dg|$.

Now we are going to use one more piece of information
available thanks to Theorem  2.1 of \cite{Kr13} which is that
$v_{K}\in W^{1,2}_{p}(\Omega_{T})$ for any $p$. Then
treating \eqref{11.1.5} near $(0,T)\times\partial \Omega$
we can flatten $\partial \Omega$ near any given point,
then continue $v-g$ (in the new coordinates)
 across the flat   boundary in an odd
way. We will then have a function of class $W^{1,2}_{d+1}$
to which Theorem 4.2.6 of \cite{Kr87} is applicable. 
In this way we   estimate the $C^{\alpha}$-norm of $v$
near the boundary of $\Omega$ and in combination with
\eqref{11.1.6} obtain that
\begin{equation} 
                                                         \label{11.1.7}
|v_{K}|_{C^{\alpha}(\bar\Omega_{T})}\leq N_{0},
\end{equation}
where $N_{0}$ depends only   on $T$, $d$, $\delta$, $K_{0}$,  the diameter
of $\Omega$, and the $L_{d+1}(\Omega_{T})$-norms of $\bar{H}$ and 
$|\partial_{t}g|+|D^{2}g|+|Dg|$.

It follows that there is a sequence $K_{n}\to\infty$ and a function
$v$ such that $v^{n}:=v_{K_{n}}\to v$ uniformly in $\bar{\Omega}_{T}$.
Of course, \eqref{11.1.2} holds,
owing to Theorem \ref{theorem 11.1.1}. Furthermore,
  \eqref{11.1.7} holds with the same
constants and $v$ in place of $v_{K}$ and $Dv^{n}\to Dv$
locally uniformly in $\Omega_{T}$.

Next, we need an analog of Lemma 6.1 of \cite{Kr13}. Introduce
$$
H_{0}(u'',t,x)=H(v(t,x),Dv(t,x),u'',t,x).
$$
\begin{lemma}
                                           \label{lemma 9.20.1} 
There is a constant $N$ depending only on $d$ and
$\delta$ such that
for any $C_{r}(t,x)$ satisfying $C_{r}(t,x)\subset\Omega_{T}$  and
$\phi\in W^{1,2}_{d+1}(C_{r}(t,x))
\cap C(\bar C_{r}(t,x))$ we have on $C_{r}(t,x)$ that
\begin{equation}
                                             \label{9.20.1}
v\leq \phi+Nr^{d/(d+1)}\|(\partial_{t}\phi+
H_{0}[\phi])^{+}\|_{L_{d+1}(C_{r}(t,x))}
+\max_{\partial'C_{r}(t,x)}(v-\phi)^{+} .
\end{equation}
\begin{equation}
                                             \label{9.20.2}
v\geq \phi-Nr^{d/(d+1)}\|(\partial_{t}\phi+
H_{0}[\phi] )^{-}\|_{L_{d+1}(C_{r}(t,x))}
-\max_{\partial'C_{r}(t,x)}(v-\phi)^{-} .
\end{equation}
\end{lemma}

Proof. Observe that
$$
-\partial_{t}\phi-\max(H_{0}[\phi] ,P[\phi]-K_{n})=
-\partial_{t}\phi-\max(H_{0}[\phi] ,P[\phi]-K_{n})
$$

$$
+\partial_{t}v^{n}+
\max(H_{0}[v^{n}] ,P[v^{n}]-K_{n})+I_{n}
$$

$$
=
\partial_{t}(v^{n}-\phi)+a_{ij}D_{ij}(v^{n}-\phi)
+I_{n},
$$
where $a=(a_{ij})$ is an $\bS_{\hat{\delta}}$-valued
function and
$$
I_{n}=\max(H [v^{n}] ,P[v^{n}]-K_{n})-\max(H_{0}[v^{n}] ,P[v^{n}]-K_{n}).
$$
Notice that  
$$
\sup_{\bar C_{r}(t,x)}|I_{n}|\leq \omega\big((\sup_{\bar C_{r}(t,x)}(|
v-v^{n}|+|Dv-Dv^{n}|)\big)\to 0
$$
as $n\to\infty$.

It follows by Theorem 3.1  of \cite{Kr78} or Theorem 3.3.9
of \cite{Kr87}
that for $r\in(0,1]$
$$
v^{n}\leq \phi+\max_{\partial'C_{r}(t,x)}(v^{n}-\phi)^{+}
$$
\begin{equation}
                                             \label{9.20.3}
+Nr^{d/(d+1)}\|(\partial_{t}\phi+I_{n}+
\max(H[\phi] ,P[\phi]-K_{n}))^{+}\|_{L_{d+1}(C_{r}(t,x))},
\end{equation}
where the constant $N=N(d,\delta)$. Actually the above references
only say that \eqref{9.20.3} holds with $N=N(r,d,\delta)$
in place of $Nr^{d/(d+1)}$.
However, the way this constant depends on $r$ is easily discovered
by using parabolic dilations.
 We obtain \eqref{9.20.1} from \eqref{9.20.3}
by letting $n\to\infty$.
In the same way \eqref{9.20.2} is established.
The lemma is proved.

After that the proof of Theorem \ref{theorem 9.23.1}
in our particular case of bounded $\bar{H}$ is achieved
in the following way.
Using \eqref{9.20.1} and repeating the proof of Theorem 2.3 
of \cite{Kr13}
(see Section 6 there), we easily obtain that,
if $(t_{0},x_{0})\in\Omega_{T}$ and $\phi\in W^{1,2}_{d+1,loc}
(\Omega_{T})$ are such that $v-\phi$ attains a local
maximum at $(t_{0},x_{0})$ and $v(t_{0},x_{0})=\phi(t_{0},x_{0})$,
then
\begin{equation}
                                                 \label{11.19.1}
 \lim _{ r\downarrow0}\esssup_{C_{r}(t_{0},x_{0}) }
\big[\partial_{t}\phi(t,x) +
H(v(t,x),Dv(t,x),D^{2}\phi(t,x),t,x)\big]\geq 0.
\end{equation}
Here $v(t,x),Dv(t,x)$ can be replaced with
$v(t_{0},x_{0}),Dv(t_{0},x_{0})$. Furthermore,
if $\phi\in W^{1,2}_{p,loc}
(\Omega_{T})$ with $p>d+2$, then by embedding theorems
$\phi\in C^{1+\alpha}_{loc}(\Omega_{T})$, where $\alpha\in(0,1)$,
 and hence
$$
(v(t_{0},x_{0}),Dv(t_{0},x_{0}))=(
\phi(t_{0},x_{0}),D\phi(t_{0},x_{0})).
$$
It follows that one can replace $v(t,x),Dv(t,x)$
with $\phi(t,x),D\phi(t,x)$
in \eqref{11.19.1} and then, by definition $v$,
is an $L_{p}$-viscosity subsolution.
 
The fact that it is also an $L_{p}$-viscosity supersolution
is proved similarly on the basis of \eqref{9.20.2}.

In  case of general $\bar H$ we introduce $u_{n}$ as the solutions
found according to Theorem  \ref{theorem 9.23.1}
of \eqref{9.23.2} in $\Omega_{T}$ with 
$$
H(u,t,x)I_{\bar{H}(t,x)\leq n}
+F(u'',t,x)I_{\bar{H}(t,x)> n}=F(u'',t,x)+G(u,t,x)I_{\bar{H}(t,x)\leq n}
$$
in place of $H(u,t,x)$ and  with the same boundary
condition $u_{n}=g$ on $\partial'\Omega_{T}$.
From the above we see that the estimates of $|u_n|_{C^{\alpha}(\bar\Omega_{T})}$
and $[u_{n}]_{C^{\kappa}(\bar C_{r}(t,x))}$ are uniform with
respect to $n$. This allows us to repeat  what was said  
about $v^{n}$ with obvious changes and brings the proof
of Theorem \ref{theorem 9.23.1} to an end.


\mysection{A minimax representation of  nonlinear
functions}
                                                       \label{section 11.4.7}

Here we complement the results of \cite{Kr11}
which originated in \cite{Kr70} by providing
a formula better suited for viewing
nonlinear PDEs as Isaacs equations.

Let $d_{1}\geq1$ be an integer.
Fix  a closed bounded subset $B$ of $\bR^{d_{1}+1}$.
Let $H(u)$ be a real-valued Lipschitz continuous
function given  on $ \bR^{d_{1}}$.
 As a Lipschitz continuous
function $H$ is differentiable on a set $D'_{H}\subset
\bR^{d_{1}} $
of full measure. 
We introduce
\begin{equation}
                                        \label{4.10.1}
\cL(H):=\{(H(u)-\langle u, D  H(u)
\rangle, D  H(u)):u\in D_{H}'\},
\end{equation}
where  $\langle\cdot,\cdot\rangle$
is the scalar product in $\bR^{d_{1}}$,
and  assume that $\cL(H)\subset B$.
Observe that for $u\in D'_{H}$ we have
$$
 H(u)-\langle u, D  H(u)
\rangle=\frac{\partial}{\partial t}[tH(u/t)]\big|_{t=1},
$$
so that, owing to the boundedness of $B$, $H$ is boundedly 
inhomogeneous.

Here is Remark 2.1 of \cite{Kr11} (in which we correct 
  an obvious misprint). 

\begin{theorem}
                              \label{theorem 4.6.1}
Under the above assumptions we have on $\bR^{d_{1}}$ that
$$
H(u)=\max_{y\in \bR^{d}}\min_{\substack{(f,l)\in B,\\
f+\langle l,y\rangle\geq H(y)}}
[f+
\langle l,u\rangle]
$$
and the sets $\{(f,l)\in B:f+\langle l,y\rangle\geq H(y)\}$
are nonempty and closed for any $y\in \bR^{d}$.

\end{theorem}

Next, let $B_{0}$ be a relatively strictly convex
closed bounded set in $\bR^{d_{1}+1}$   such that $B_{0}\supset B$
and the distance between the relative boundaries
of $B$ and $B_{0}$ is strictly positive. 
Then introduce $A:=\bR^{d_{1}}$ and for
 $\alpha\in A$   define
$$
G(\alpha)=\sup_{(f,l)\in B_{0}}(f+\langle l,\alpha\rangle).
$$
Next,
let $\cP$ be the smallest hyperplane containing $B_{0}$,
and, by using the assumption about the boundaries
of $B$ and $B_{0}$,
define $\Gamma$ as a closed (in the topology of $\cP$)
 convex subset of $\cP$
with the origin lying in the relative 
(in the topology of $\cP$)
interior of $\Gamma$
such that $(f,l)+\Gamma\subset B_{0}$ for any 
$(f,l)\in B$. Define
$$
\gamma(u)=\sup_{(f,l)\in \Gamma}[f+\langle l,u\rangle]
$$
and observe that since $\mu(\pm B_{0})\subset\Gamma$
for a constant $\mu>0$ we have that
$$
\mu|G(\alpha)|\leq\gamma(\alpha),\quad
\mu|H(\alpha)|\leq\gamma(\alpha).
$$
Furthermore, since $(f,l)+\Gamma\subset B_{0}$ for any $(f,l) \in B $
we have that for any $(f,l) \in B $
\begin{equation}
                                              \label{9.26.1}
f+\langle l,\alpha\rangle+\gamma(\alpha)\leq
\sup_{(f',l')\in B_{0}}(f'+\langle l',\alpha\rangle),
\quad G(\alpha)
-[f+\langle l,\alpha\rangle]\geq \gamma(\alpha),
\end{equation}
which shows that
$$
\lambda^{\alpha\beta}_{H}:=1\wedge\frac{G(\alpha)-H(\alpha)}{G(\alpha)
-[f+\langle l,\alpha\rangle]} 
\quad\quad\bigg(\frac{0}{0}=1\bigg)
$$
 is well defined for  
 $\alpha\in A$ and $\beta=(f,l)\in  B$ 
and, of course, $\lambda^{\alpha\beta}_{H}\in[0,1]$.

Next, observe that the graph of
 $$
G(\xi,\alpha):=\sup_{(f,l)\in B_{0}}(f\xi
+\langle l,\alpha\rangle)
$$
is a cone with respect to $(\xi,\alpha)$ which
is once continuously differentiable with respect to
$(\lambda,\alpha)$ everywhere apart from the origin
due to the strict convexity of $B_{0}$. Since the plane
$\xi=1$ does not pass through the origin,
$G(\alpha)$ is once continuously differentiable.

Now for $\alpha\in A$ and $\beta=(f,l)\in  B$ set
$$
f^{\alpha\beta}_{H}=
\lambda^{\alpha\beta}_{H}f+(1-\lambda^{\alpha\beta}_{H})[
G(\alpha)-\langle\alpha,DG(\alpha)\rangle],
$$
$$
l^{\alpha\beta}_{H}=\lambda^{\alpha\beta}_{H}l+
(1-\lambda^{\alpha\beta}_{H})DG(\alpha).
$$

Obviously $(f^{\alpha\beta}_{H},l^{\alpha\beta}_{H})\in B_{0}$.
In this way on the set $\cH(B)$ of functions $H$
satisfying the assumptions stated in the beginning of the section
 we constructed a mapping   sending each
$H\in\cH(B)$ into the function $(f^{\alpha\beta}_{H}
,l^{\alpha\beta}_{H})$ defined  on $A\times B$

\begin{theorem}
                                           \label{theorem 8.30.1}
For any $H\in\cH(B)$ and any $u\in\bR^{d_{1}}$ we have
$$
H(u)=\supinf_{\alpha\in A\,\,\beta\in B}
[f^{\alpha\beta}_{H} +
\langle l^{\alpha\beta}_{H} ,u\rangle].
$$
Furthermore, if $H,F\in\cH(B)$, then for any $\alpha\in A$
and $\beta=(f,l)\in B$ we have
$$
|f^{\alpha\beta}_{H}-f^{\alpha\beta}_{F}|\leq
\frac{|H(\alpha)-F(\alpha)|}{\gamma(\alpha)}|f+\langle \alpha,
DG(\alpha)\rangle-G(\alpha)|\quad\quad\bigg(\frac{0}{0}=0\bigg),
$$
$$
|l^{\alpha\beta}_{H}-l^{\alpha\beta}_{F}|\leq
\frac{|H(\alpha)-F(\alpha)|}{\gamma(\alpha)}|l-
DG(\alpha)\rangle |\quad\quad\bigg(\frac{0}{0}=0\bigg).
$$
\end{theorem}

Proof.  
Observe that, for
$\beta=(f,l)\in B$,  if $f+\langle l,\alpha\rangle
\geq H(\alpha)$, then
 $G(\alpha)-H(\alpha)
\geq G(\alpha)-[f+\langle l,\alpha\rangle]$ 
and  $\lambda_{H}^{\alpha\beta}=1$ (no matter $\gamma(\alpha)=0$
or $\gamma(\alpha)>0$) and $f=f^{\alpha\beta}_{H}$
and $l=l^{\alpha\beta}_{H}$.
It follows that
$$
\min_{\substack{\beta=(f,l)\in B,\\
f+\langle l,\alpha\rangle\geq H(\alpha)}}[f+
\langle l,u\rangle]\geq
\inf_{\substack{\beta\in B,\\
\lambda^{\alpha\beta}_{H}=1}}[f^{\alpha\beta}_{H}+
\langle l^{\alpha\beta}_{H},u\rangle]
\geq
\inf_{\beta\in B}(
 f^{\alpha\beta}_{H}+\langle l^{\alpha\beta}_{H},u\rangle).
$$

Furthermore, for
$\beta=(f,l)\in B$,  if $\lambda^{\alpha\beta}_{H}=1$,
then  $G(\alpha)-H(\alpha)\geq G(\alpha)
-[f+\langle l,\alpha\rangle]$, so that
$(f^{\alpha\beta}_{H},l^{\alpha\beta}_{H})=(f,l)\in B_{0}$ and 
$$
f^{\alpha\beta}_{H}+\langle l^{\alpha\beta}_{H},
\alpha\rangle\geq H(\alpha).
$$
In addition, for
$\beta=(f,l)\in B$,  if   $\lambda_{H}^{\alpha\beta}<1$, then
as always
$(f^{\alpha\beta}_{H},l^{\alpha\beta}_{H}) \in B_{0}$ and
$$
 f^{\alpha\beta}_{H}+\langle l^{\alpha\beta}_{H},\alpha\rangle
=\lambda^{\alpha\beta}_{H}[f+\langle l,\alpha\rangle]
+(1-\lambda^{\alpha\beta}_{H})G(\alpha)=H(\alpha).
$$
Hence,
$$
\min_{\substack{\beta=(f,l)\in B,\\
f+\langle l,\alpha\rangle\geq H(\alpha)}}[f+
\langle l,u\rangle]\geq
\inf_{\beta\in B}(
 f^{\alpha\beta}_{H}+\langle l^{\alpha\beta}_{H},u\rangle)
\geq
\min_{\substack{ (f,l)\in B_{0},\\
f+\langle l,\alpha\rangle\geq H(\alpha)}}[f+
\langle l,u\rangle]
$$
and the first assertion  of the theorem follows from
Theorem \ref{theorem 4.6.1}. 

To prove the second assertion it suffices to note that,
for instance.
$$
f^{\alpha\beta}_{H}-f^{\alpha\beta}_{F}=
(\lambda^{\alpha\beta}_{H}-\lambda^{\alpha\beta}_{F})
[f+\langle \alpha,
DG(\alpha)\rangle-G(\alpha)],
$$
where the right-hand side is zero if $\gamma(\alpha)=0$,
and after that use \eqref{9.26.1}.
The theorem is proved.

\begin{theorem}
                                   \label{theorem 9.4.1}
Let $H$ also depend on  parameters $(t,x)\in\bR^{d+1}$
and let it satisfy the assumptions in the beginning
of the section for each $(t,x)$. Assume that
$H(u,t,x)$ is measurable with respect to $(t,x)$ and
there is a function $\bar{H}(u,t)$ also satisfying
the assumptions in the beginning
of the section for each $t$ and  measurable with respect to $t$.
 Denote
$$
\theta=\int_{C_{1,1}}\sup_{u:\gamma(u)\ne0}
\frac{|H(u,t,x)
-\bar{H}(u,t)|}{\gamma(u)}\,dxdt.
$$
Also let $(f^{\alpha\beta}(t,x),l^{\alpha\beta}(t,x))$
correspond to $H(u,t,x)$ and $(\bar{f}^{\alpha\beta}(t),
\bar{l}^{\alpha\beta}(t))$
correspond to $\bar{H}(u,t)$ constructed as  before
Theorem \ref{theorem 8.30.1}. Then
$$
\int_{C_{1,1}}\sup_{\alpha\in A,\beta\in B}
\big(|f^{\alpha\beta}(t,x)-\bar{f}^{\alpha\beta}(t)|
+|l^{\alpha\beta}(t,x)-
\bar{l}^{\alpha\beta}(t)|\big)\,dxdt\leq N\theta,
$$
where the constant $N$ depends only on $B$ and $B_{0}$.

\end{theorem} 

Proof. Let $\lambda^{\alpha\beta}(t,x)$
correspond to $H(u,t,x)$ and $\bar\lambda^{\alpha\beta}(t)$
correspond to $\bar{H}(u,t)$ constructed as  before
Theorem \ref{theorem 8.30.1}. Then it suffices to prove that
$$
\int_{C_{1,1}}\sup_{\alpha\in A,\beta\in B}|
\lambda^{\alpha\beta}(t,x)-\bar\lambda^{\alpha\beta}(t)|\,dxdt
\leq N\theta.
$$
For $\beta=(f,l)\in B$ we have
$$
|\lambda^{\alpha\beta}(t,x)-\bar\lambda^{\alpha\beta}(x)|\leq
\bigg|\frac{G(\alpha)-H(\alpha,t,x)}{G(\alpha)
-[f+\langle l,\alpha\rangle]}-
\frac{G(\alpha)-\bar{H}(\alpha,t)}{G(\alpha)
-[f+\langle l,\alpha\rangle]}\bigg|
$$
$$
\leq \frac{|H(\alpha,t,x)
-\bar{H}(\alpha,t)|}{\gamma(\alpha)}
$$
(with $0/0=0$)
and our assertion follows. The theorem is proved.

\end{document}